\documentclass[11pt]{article}
\usepackage{amsmath, amsthm, amssymb}    % May not all be necessary
\usepackage{bridges}                     % Custom bridges proceedings style
\usepackage{graphicx}                    % For including pictures
\usepackage{hyperref}                    % For formatting (clickable) URLs
\usepackage[english]{babel}
\usepackage{multicol}
\usepackage{blindtext}

\usepackage[latin1]{inputenc}
\usepackage{enumerate}
\usepackage{psfrag}

\theoremstyle{plain}
\newtheorem{theorem}{Theorem}%[section]

\theoremstyle{definition}

\numberwithin{equation}{section}

\newcommand{\Z}{\mathbb{Z}}

\newcommand{\abs}[1]{\left| #1 \right|}

\title{Toss and Spin Juggling State Graphs}
\author{Harri Varpanen\\
Dept.~of Mathematics and Systems Analysis, Aalto University\\
P.O.Box 11100, FI-00076 Aalto, Finland\\
{\tt harri.varpanen@aalto.fi}
}
\date{}				% Suppress any date on submissions

\begin{document}

\maketitle

% Prevent page number 1 from being printed on the first page.
\thispagestyle{empty}

\begin{abstract}
We review the state approach to toss juggling and extend the approach
to spin juggling, a new concept. We give connections to current research
on random juggling and describe a professional-level juggling performance that
further demonstrates the state graphs and their research.
\end{abstract}

\section*{Introduction}
A well-studied mathematical model for juggling called {\it siteswap} 
has been described in previous Bridges papers by Bracken \cite{Bracken_2008}
and Naylor \cite{Naylor_2012}. The model is essentially
one-handed and it only keeps track of the throwing schedules of
the objects thrown one by one. While the model can be generalized to include
several hands and simultaneous throws, it lacks any spatial information.

Juggling is commonly viewed as {\em tossing} various objects: Balls, rings, clubs, and their
variants such as torches and chainsaws. The siteswap model is purely temporal and regards all these objects as {\em particles}, and accordingly the mathematical treatment
is mostly combinatorial in nature. In practice spatial information is often
included separately, for example: ''Siteswap pattern
$423$ with clubs; 4s as double spins, 2s as swings, and 3s as single spins
thrown behind the back.''

However, there are other common forms of juggling that are more
spatial in nature and deserve a
mathematical treatment of their own.
Although the word spin can be used as a club rotation count like above,
we introduce {\em spin juggling} to include poi spinning,
devilstick manipulation, and diabolo juggling.
Tossing may occur in these
forms of juggling as well, but it is less significant
compared to the spinning that consists of various spatial
configurations and rotational transitions between them. 

In this paper we review the state approach to toss juggling,
and introduce a {\em poi state graph} as an example of a {\em spin
layer} that can be added to toss juggling. We describe
the mathematical research related to random walks on juggling state graphs
and, finally, depict a juggling performance built on the
material in this paper.

\section*{Historical background}
The mathematical modeling of juggling patterns started around 1980,
when computers gained popularity and scholars envisioned robots
and simulators. Shannon \cite{Shannon_1980} is
probably the first mathematical manuscript on juggling, with emphasis
on robotics.

The siteswap model was independently invented by several people in the early
1980s and was first published by Magnusson and Tiemann
\cite{Magnusson_Tiemann_1989}. Since then a dozen or so serious papers,
in addition to a monograph by Polster \cite{Polster_2003},
have been published on the subject. We mention the
following papers:

\newpage 

\begin{itemize}
\item Graham et al. \cite{Buhler_et_al_1994, Butler_Graham_2010, Chung_et_al_2010, Chung_Graham_2008} on the combinatorics of periodic siteswaps
\item Ehrenborg and Readdy \cite{Ehrenborg_Readdy_1996} pointing
out a connection between $n$-periodic siteswaps and the affine
Weyl group $\widetilde{A}_{n-1}$ (along with further interesting combinatorics)
\item Warrington \cite{Warrington_2005} on uniformly random siteswaps with
bounded throws
\item Devadoss and Mugno \cite{Devadoss_Mugno_2007} on generating braids via siteswaps
\item Leskel\"a and Varpanen \cite{Leskela_Varpanen_2012}
on general random siteswaps
\item Knutson, Lam and Speyer \cite{Knutson_Lam_Speyer_2011} decomposing
the Grassmannian manifold of $k$-planes in $n$-space via $n$-periodic
siteswaps with $k$ particles.
\end{itemize}
There are also various toss juggling simulators such as Juggling Lab
\cite{JugglingLab} coordinated by Jack Boyce.

\section*{Siteswap state graphs}
Siteswap juggling patterns are usually modeled as bijections
$f \colon \Z \to \Z$ with $f(t) \ge t$ for each $t \in \Z$, the
interpretation being that a juggling pattern is eternal, a particle thrown at
beat $t$ is scheduled to be next thrown at a future beat $f(t)$, and that
at most one particle is allowed to be thrown at a given beat.
Here the number
$f(t) - t \ge 0$ is called the {\it throw} at beat $t$, and a zero-throw
is interpreted as an empty throw where the juggler has nothing to throw and
waits. A nontrivial orbit of $f$ corresponds to the throw times
of a particle, and the long-time average of the throws equals the number of
particles. The model is
essentially one-handed, but it also comprises the commonly encountered
patterns with alternating hands.

Figure 1 depicts a three-particle, period-three pattern.
Instead of beat numbers $t$ we have listed the {\it throw values} $f(t)-t$; the pattern is accordingly called $450$. 

\begin{figure}[h!tbp]
\begin{center}
\includegraphics[scale=1.7]{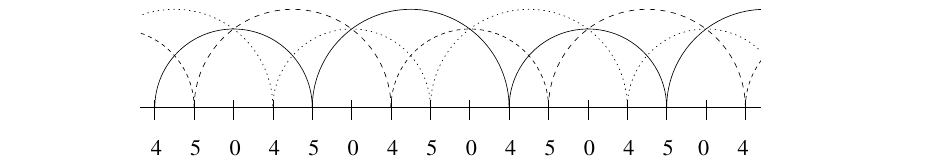} \\
\label{state1}
\end{center}
\caption{The pattern $450$.}
\end{figure}

\noindent Just before a $5$-throw is performed in Figure 1, the three
particles in the air are
scheduled to be next thrown zero, two and three beats from ''now'', the
vertical line in Figure 2.
We say that the juggler is at {\it state} $\{0,2,3\}$. 

\begin{figure}[h!tbp]
\begin{center}
\includegraphics[scale=1.4]{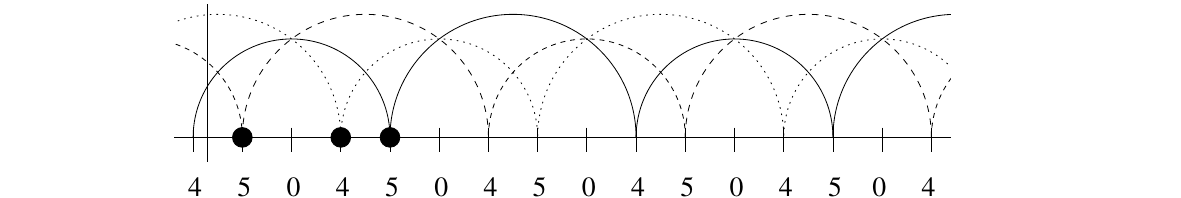} \\
\end{center}
\label{state2}
\caption{A state $\{0,2,3\}$ in the pattern $450$.}
\end{figure}

As the $5$-throw is performed and time advances one beat forward, the state
changes from $\{0,2,3\}$ to $\{1,2,4\}$ (Figure 3):
\newpage
\begin{figure}[h!tbp]
\begin{center}
\includegraphics[scale=1.4]{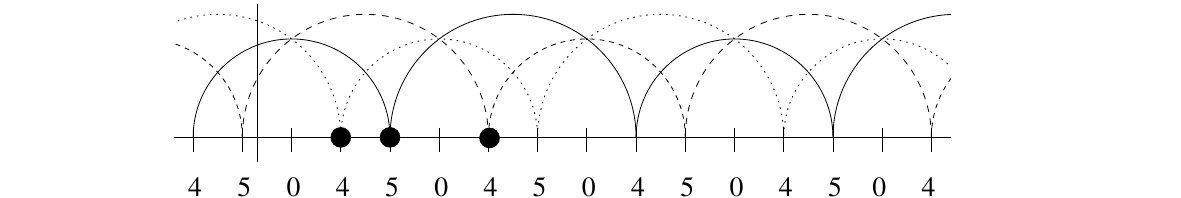} \\
\end{center}
\label{state3}
\caption{A state $\{1,2,4\}$ in the pattern $450$.}
\end{figure}
The pattern $450$ can hence be represented via its {\it state cycle}
as in Figure 4 below. The states
have been drawn vertically to reflect particles falling down and thrown
up. The states $\{0,1,3\}$, $\{0,2,3\}$ and
$\{1,2,4\}$ are three-element subsets of the five-element set $\{0,1,2,3,4\}$,
because there are three particles and the maximum throw value is five.
The throw values $4$, $5$ and $0$ yield transitions from one state to the next.
\begin{figure}[h!tbp]
\begin{center}
\includegraphics[scale=0.9]{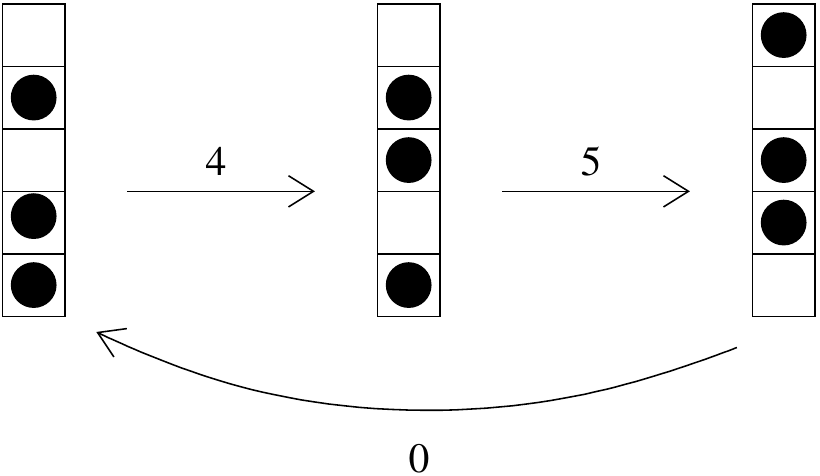} \\
\end{center}
\label{450_cycle}
\caption{The state cycle for the pattern $450$.}
\end{figure}

Now any three-particle pattern with maximum throw five is a path in the following {\it state graph} (Figure 5):

\begin{figure}[h!tbp]
\begin{center}
\includegraphics[scale=1.3]{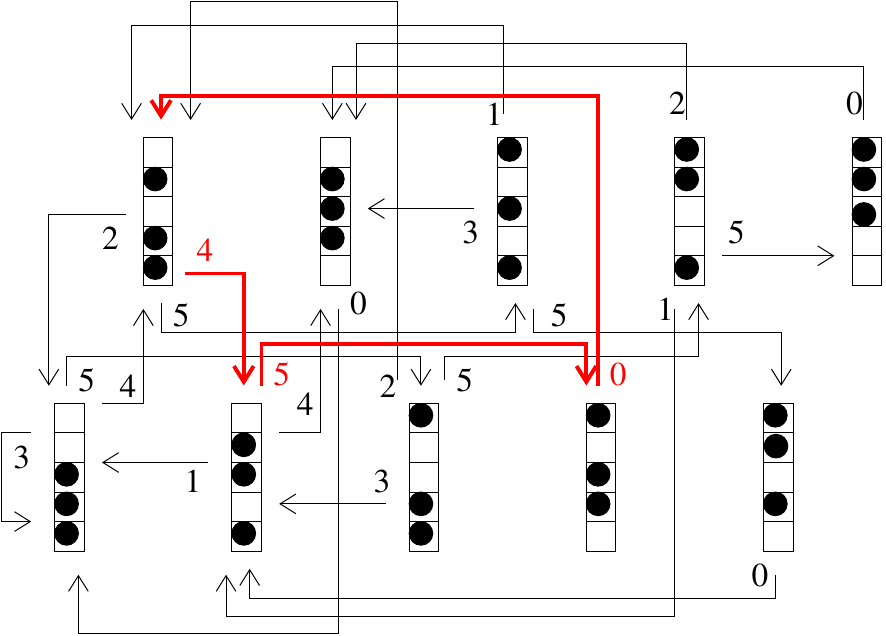} \\
\end{center}
\label{stategraph1}
\caption{The state graph for three particles and maximum throw five. The pattern $450$ is highlighted.}
\end{figure}

The state model was independently invented by Jack Boyce and Allen Knutson
in 1988 and is widely used by jugglers to communicate juggling patterns
and transitions between them.
The model naturally generalizes to allow
simultaneous throws ({\em multiplex juggling}) and many hands ({\em multihand
juggling}) as explained in Polster \cite{Polster_2003}.
%as in Figure 6, where particles fall two
%levels per beat for consistency with the alternating-hands case, and where
%x denotes a crossing throw.
%
%\begin{figure}[h!tbp]
% \begin{center}
%\hspace{4ex} \includegraphics[scale=0.9]{slam_multi04.pdf} 
%\end{center}
%\label{multi1}
%\caption{A five-particle pattern with two hands throwing simultaneously.
%The reader is invited to enter the pattern $([4x2],[22])([22],[4x2])([6x4],[42])*$ into the Juggling Lab simulator.}
%\end{figure}

\section*{Spin juggling}

\begin{multicols}{2}

The simplest form of spin juggling is {\em poi spinning}, where a juggler
spins tethered weights, one in each hand, through a variety of rhythmical
and geometric patterns. Figure 6 shows a schematic poi juggler in
its simplest state: The tethers are not braided.
The weights rotate 180 degrees per beat in paths that resemble half-circles:
One weight from down to up, the other from up to down. Braiding of
the tethers occurs when one of the weights moves from one side of the
juggler's body to the other.

\columnbreak

%\begin{figure}[h!tpb]
\begin{center}
\includegraphics[scale=0.8]{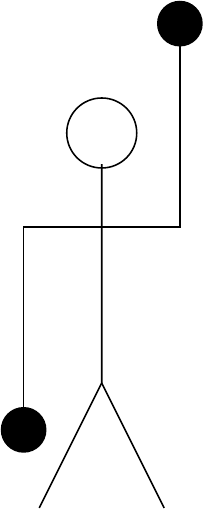} \\
{\small {\bf Figure 6} : {\em A poi juggler.} }
\end{center}
%\caption{test}
%\end{figure}
\stepcounter{figure}

\end{multicols}

In order to capture the essence of poi spinning, we ignore the weights
and regard the tethers as extended hands (made of a rubber string, say).
We also ignore the juggler's body except from the fixed shoulder line.
We assume that the juggler is viewed from behind and that the direction
of the spinning is forward: The rotation of the weight from up to down
takes place in front of the juggler's body, rather than behind. We also
make the simplifying assuption that a hand may move to the other side of the
body (shoulder line) only while rotating from down to up.

\begin{figure}[h!tbp]
\begin{center}
\includegraphics[scale=0.5]{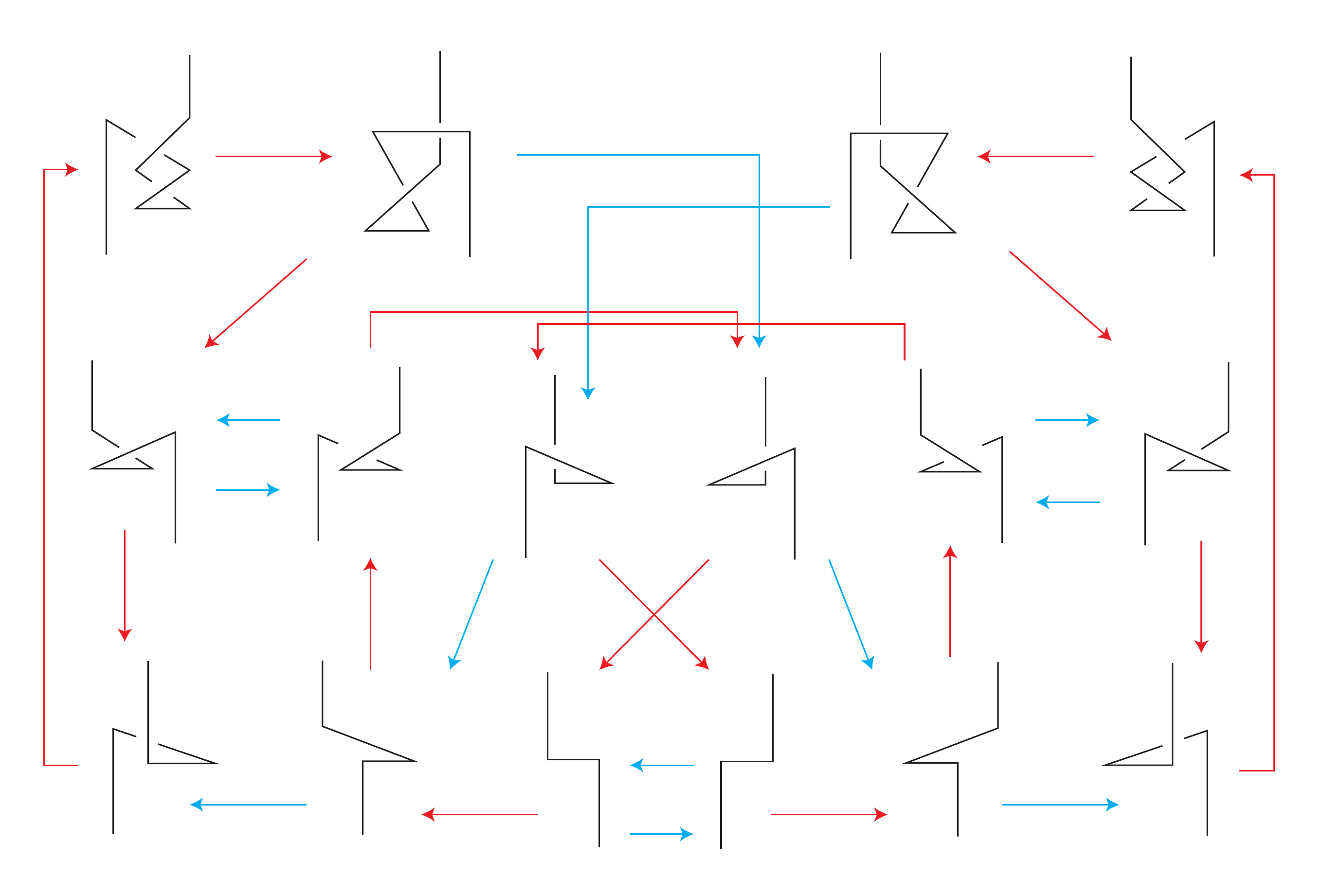} 
\end{center}
\caption{A poi state graph.}
\end{figure}

\newpage

The resulting poi state graph is given in Figure 7, where each state
(black line) consists of a fixed horizontal shoulder line along with a
changing hand configuration of the juggler.
The arrows indicate the 180-degree rotations of the hands from down to up and
vice versa. A red arrow indicates that a hand moves to the other side of
the body, while a blue arrow indicates that no such move takes place. (Remember
that only the hand rotating from down to up is allowed to move to the other
side of the body.) The braiding of the hands is bounded to a maximum of two crossings, enough to cover the commonly encountered ''3-weave'' and ''5-weave''
poi spinning patterns.

For example, if a move is denoted by R (red) and a non-move by B (blue),
then the ground pattern formed by the two states in the middle bottom row of
Figure 7 is
$\ldots$BBB$\ldots$ or, assuming that the symbols repeat
periodically, simply B.
The 3-weave pattern is similarly RRB and
% \mbox{$\ldots$RRBRRBRRB$\ldots$} or simply RRB, and the 5-weave
the 5-weave pattern is BRRBB, but these patterns can not be started directly from
the ground pattern. Starting from the
% $\ldots$BRRBBBRRBBBRRBB$\ldots$ or simply BRRBB.
ground pattern and including the default transition moves in parentheses,
the 3- and 5-weave patterns read (R)RRB(RRR) and (R)BRRBB(BRRBR), respectively.

\begin{multicols}{2}

\begin{center}
\includegraphics[scale=0.2]{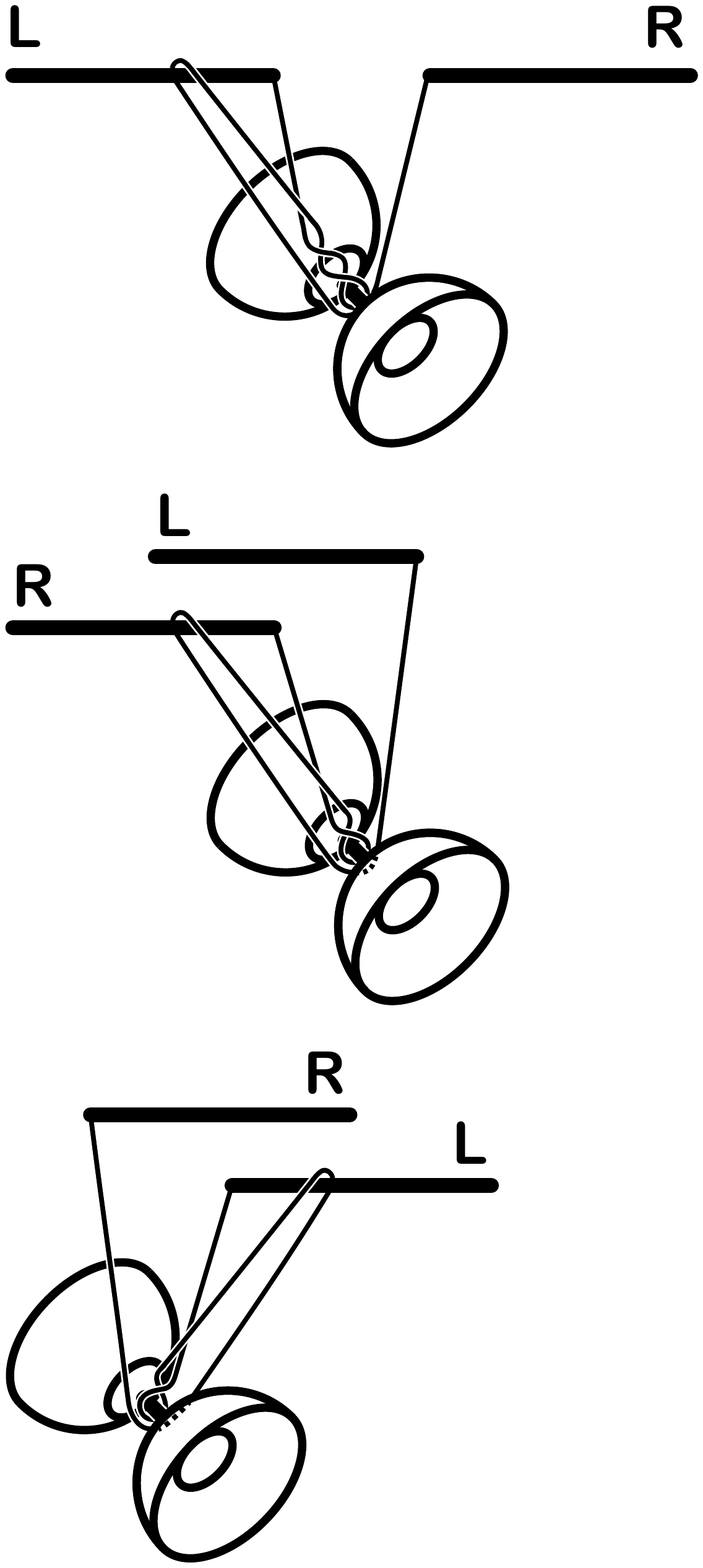} \\
{\small {\bf Figure 8} : {\em A few diabolo states.} }
\end{center}

\columnbreak

The poi hand movements form a spatial layer that can be combined with siteswap.
For example, the siteswap 522 combined with the 3-weave RRB is a nice
pattern when juggled with clubs. In addition, the poi hand movements
are similar to the hand movements that occur in devilstick juggling or
staff manipulation (these are considered mathematically equivalent), so the
poi state graph applies to these forms of juggling as well. Finally,
similar state graphs have been considered
for the diabolo by M\"annist\"o \cite{Mannisto_2012}. Figure 8 shows
some states for the diabolo, the complete state graph being too complicated
to be explained in this paper. 

%\begin{figure}[h!tpb]
%\caption{test}
%\end{figure}
\stepcounter{figure}

\end{multicols}

%\begin{multicols}{2}
%\begin{figure}[h!tbp]
%\begin{center}
%\includegraphics[scale=0.2]{Samuli_kuvat2.pdf} 
%\end{center}
%\end{figure}
%\columnbreak 
%\mbox{ } 
%Figure X. A few diabolo states. The full state graph is too complicated to be shown here.
%\end{multicols}

\section*{Random juggling}
Although no one is able to perform truly random juggling, one may still
consider random juggling from a mathematical point of view.
Random walks on juggling state graphs are fascinating because
a juggling state graph possesses a clear mathematical structure:
Nontrivial, but not too complicated to produce exact
results.

Given a juggling state, consider all its admissible follower states. For
example, the leftmost state in Figure 5 has three admissible followers
(reached by transitions labeled 3, 4, and 5), one of which is itself.
Associate to each follower state a probability such that the
probabilities sum to one. In other words, choose a transition
from the current state to the next state by throwing a dice. 
Assume that the dice is known at each state; then ask:
''{\em How popular, relative to
each other, are the juggling states in the long run?}''
In other words, one asks for the {\em steady-state distribution} (or {\it equilibrium distribution}) of the random juggling process.

While it is often possible to calculate numerical approximations of equilibrium
distributions, explicit formulas are typically considered difficult or
impossible. However, the structure present in juggling state graphs allows
for some explicit formulas, the first of which was obtained by Warrington
\cite{Warrington_2005}
in the case where each siteswap throw is chosen uniformly randomly:
If there are, say, three alternatives for the throw whenever the juggler does not wait (as is the case in Figure 5), then each alternative is chosen with probability one
third. The exact long-term visit frequencies are shown in Figure 9 and
calculated using the formula in Theorem 1.

\begin{theorem}[\cite{Warrington_2005} Thm. 5]
In uniform random juggling with $k$ particles and maximum throw $m$, the frequency 
 of a state $B \subset \{0,1,\ldots,m-1\}$, $\abs{B} = k$ is
\[
 \left\{\begin{smallmatrix}m+1\\m+1-k\end{smallmatrix}\right\}^{-1}\prod_{x \in B} | \ \{x,\ldots,m\} \setminus B \ |, 
\]
where $\left\{\begin{smallmatrix}m+1\\m+1-k\end{smallmatrix}\right\}$ is a Stirling number of the second kind.
\end{theorem}

\begin{figure}[h!tbp]
\begin{center}
\includegraphics[scale=1.2]{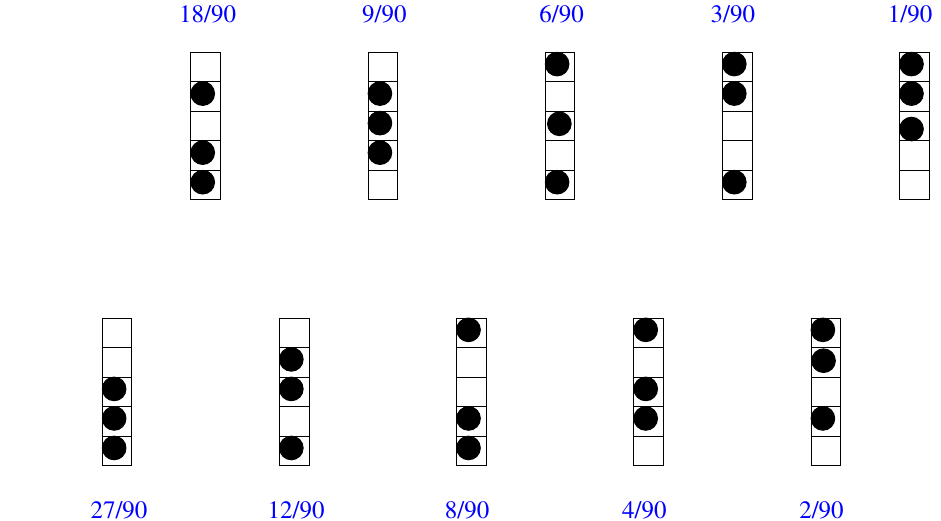}
\end{center}
\label{steady1}
\caption{The visit frequencies of the states in uniformly random siteswap
juggling with three particles and maximum throw five. The (leftmost) state
with lowest potential energy is visited almost every third beat whereas the
(rightmost) state with highest potential energy is visited only once in every 90 beats.}
\end{figure}

In our recent research \cite{Leskela_Varpanen_2012, Engstrom_Leskela_Varpanen_2013} we have considered situations where the dice is not fair and where the
throws can be arbitrary high. Explicit formulas are available in these
situations as well, and new kinds of combinatorial limit formulas emerge
when the maximum throw approaches infinity. Our ongoing research involves
juggling with antimatter and mean field juggling involving spinning hands.
These are further demonstrated in my presentation.

\section*{The presentation}

I started studying mathematics in 1994 and found juggling in 1998. Since
then my main interest has been in the mathematical models of both toss
and spin juggling. I have kept practicing regularly, and juggle fluently
(but not simultaneously) with the six most essential toss and
spin equipment: balls, rings, clubs, poi, devilstick, and diabolo.

I have built a half-hour mathematical juggling performance around the topics 
surveyed in this paper. The performance has the form of a lecture, but the
speech is prerecorded in the slides that I change with a foot pedal while
quietly demonstrating various juggling patterns; see the screenshot in
Figure 10. A three-minute trailer of the performance in Figure 10 is available
(as of \today) at
\begin{center}
\url{http://www.youtube.com/watch?v=u2mErXtXqMc}
\end{center}
(or search ''Varpanen Heureka 2013''). The trailer is in Finnish with English
subtitles, but the performance is available in English as well.

My aim is to illustrate for a general scientific audience in an
artistic, entertaining and graspable way how hands-on juggling models relate
to advanced mathematical research.
During the last ten years I have constantly improved the show and performed
it on various occasions such as college mathematics camps, science slam
events, mathematics teachers' conferences, science center mathematics days,
and university mathematics course lectures.
The performance can be carried out in any normal auditorium or lecture hall
with a high enough ceiling and some free floor space. 

\begin{figure}[h!tbp]
        \begin{center}
        \includegraphics[width=0.8\textwidth]{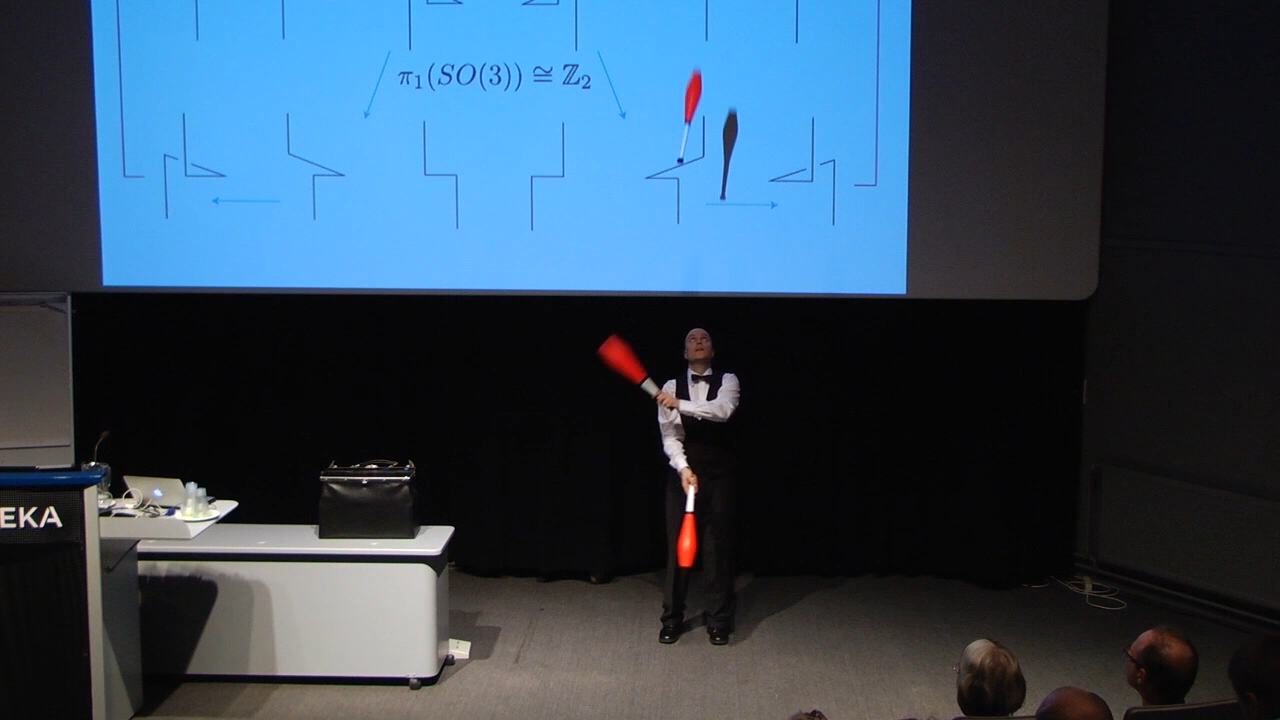}
        \end{center}
        \caption{A screenshot of my juggling lecture from 2013. The white
foot pedal can be seen on the floor at the desk.}
        %\label{FigSample}
\end{figure}

\vspace{2ex}

\section*{Conclusion}
State graphs (also called state diagrams) occur everywhere: Population
dynamics, the Internet, quantum physics, factories, etc. Juggling state graphs
are free from industrial restrictions and give rise to rich, deep and beautiful
mathematics. In addition, juggling provides an exceptional means of
communicating the mathematics to a wider audience.

While toss juggling state graphs have been analyzed for some decades,
spin juggling is a new concept that adds spatial and rotational dimensions
to juggling. Currently its mathematical formulation is still in its infancy,
and no serious spin juggling simulators exist (to the best of my knowledge).
Besides mathematical papers on spin juggling, I would love to witness a
scientific juggling simulator that generates and analyzes random walks on
both toss and spin juggling state graphs, while showing the patterns on
the screen.

\section*{Acknowledgements}
I am indebted to Alexander Engstr\"om and Lasse Leskel\"a for
mathematical support, \mbox{Joni Tikkanen} for building
the foot pedal, and Samuli M\"annist\"o for his interest
in spin juggling state graphs. I am also grateful to the anonymous Bridges
referees that largely helped in improving this paper.


\newpage

\begin{thebibliography}{99}

\bibitem{Bracken_2008} C. Bracken.
\newblock Counting The Number of Site Swap Juggling Patterns with Respect to Particular Ceilings.
\newblock {\em Proceedings of Bridges 2008: Mathematics, Music, Art, Architecture, Culture}, 219--224, 2008.

\bibitem{Buhler_et_al_1994} J. Buhler, D. Eisenbud, R. Graham and C. Wright.
\newblock Juggling drops and descents.
\newblock {\em Amer. Math. Monthly} 101~(6), 507--519, 1994.

\bibitem{Butler_Graham_2010}
S. Butler and R. Graham.
\newblock Enumerating (multiplex) juggling sequences.
\newblock {\em Ann. Comb.} 13~(4), 413--424, 2010.

\bibitem{Chung_et_al_2010}
F. Chung, A. Claesson, M. Dukes and R. Graham.
\newblock Descent polynomials for permutations with bounded drop size.
\newblock {\em European J. Combin.} 31~(7), 1853--1867, 2010.

\bibitem{Chung_Graham_2008}
F. Chung and R. Graham.
\newblock Primitive juggling sequences.
\newblock {\em Amer. Math. Monthly} 115~(3), 185--194, 2008.

\bibitem{Devadoss_Mugno_2007}
S. Devadoss and J. Mugno.
\newblock Juggling braids and links.
\newblock {\em Math. Intelligencer} 29~(3), 15--22, 2007.

\bibitem{Ehrenborg_Readdy_1996} R. Ehrenborg and M. Readdy.
\newblock Juggling and applications to $q$-analogues.
\newblock {\em Discrete Math.} 157, 107--125, 1996.

\bibitem{Engstrom_Leskela_Varpanen_2013} A. Engstr\"om, L. Leskel\"a and H. Varpanen.
\newblock Geometric juggling with $q$-analogues. 
\newblock Submitted. Preprint available at {\em arXiv:1310.2725.}

%\bibitem{Engstrom_Leskela_Varpanen_2014} A. Engstr\"om, L. Leskel\"a and H. Varpanen.
%\newblock Virtual juggling recurrence.
%\newblock In preparation.

%\bibitem{Feynman_1949}
%R. P. Feynman.
%\newblock The Theory of Positrons.
%\newblock {\em Phys. Rev.} 76, 749--759, 1949.
%
%\bibitem{Klimenko_2014}
%A. Y. Klimenko and U. Maas.
%\newblock One antimatter -- two possible thermodynamics.
%\newblock {\em Entropy} 16, 1191--1210, 2014, to appear. doi:10.3390/e16031191

\bibitem{JugglingLab}
\newblock The Juggling Lab simulator.
\newblock Online at \url{http://jugglinglab.sourceforge.net} (as of \today).

\bibitem{Knutson_Lam_Speyer_2011}
A. Knutson, T. Lam and D. E Speyer.
\newblock Positroid varieties: Juggling and geometry.
\newblock {\em Compos. Math.} 149~(10), 1710--1752, 2013.

\bibitem{Leskela_Varpanen_2012} L. Leskel\"a and H. Varpanen.
\newblock Juggler's exclusion process.
\newblock {\em J. Appl. Probab.} 49~(1), 266--279, 2012.

\bibitem{Magnusson_Tiemann_1989} B. Magnusson and B. Tiemann.
\newblock The physics of juggling.
\newblock {\em Phys. Teach.} 27, 584--588, 1989.

\bibitem{Mannisto_2012} S. M\"annist\"o.
\newblock {\em A graphical juggling notation system.}
\newblock Master's Thesis (in Finnish), University of the Arts Helsinki, 2012.

\bibitem{Naylor_2012} M. Naylor.
\newblock The mathematical art of juggling: using mathematics to predict,
describe and create.
\newblock {\em Proceedings of Bridges 2012: Mathematics, Music, Art, Architecture, Culture}, 33--40, 2012.

\bibitem{Polster_2003}
B. Polster.
\newblock {\em The Mathematics of Juggling.}
\newblock Springer Verlag, 2003.

\bibitem{Shannon_1980}
C. Shannon.
\newblock Scientific aspects of juggling.
\newblock Manuscript from ca. 1980. Published in {\em Claude Elwood Shannon:
Collected Papers}, pp. 850--864, Wiley 1993.

\bibitem{Warrington_2005}
G. S. Warrington.
\newblock Juggling probabilities.
\newblock {\em Amer. Math. Monthly} 112~(2), 105--118, 2005.

\end{thebibliography}
\end{document}